\documentclass[11pt]{article}
\usepackage{amsmath}
\usepackage{amssymb}
\usepackage{amscd}
\usepackage{xspace}
\usepackage{verbatim}
\usepackage{latexsym}
\newcommand{\mysection}[1]{
\section{#1}\setcounter{equation}{0}}
\title{\bf On the equation $-\Delta u+e^{u}-1=0$ with measures as  boundary data}
\author{
 {\bf Laurent V\'eron}\\[2mm]
{\small Laboratoire de Math\'ematiques et Physique Th\'eorique, }\\
{\small  Universit\'e Fran\c{c}ois Rabelais,  Tours,  FRANCE}}

\date{}
\begin{document}
\maketitle
{\small {\bf Abstract} If  $\Omega$ is a  bounded domain in $\mathbb R^N$, we study conditions on a Radon measure $\mu$ on $\partial\Omega$ for solving the equation $-\Delta u+e^{u}-1=0$ in $\Omega$ with $u=\mu$ on $\partial\Omega$. The conditions are expressed in terms of Orlicz capacities.
}

\vspace{1mm}
\hspace{.05in}

\newcommand{\txt}[1]{\;\text{ #1 }\;}
\newcommand{\tbf}{\textbf}
\newcommand{\tit}{\textit}
\newcommand{\tsc}{\textsc}
\newcommand{\trm}{\textrm}
\newcommand{\mbf}{\mathbf}
\newcommand{\mrm}{\mathrm}
\newcommand{\bsym}{\boldsymbol}
\newcommand{\scs}{\scriptstyle}
\newcommand{\sss}{\scriptscriptstyle}
\newcommand{\txts}{\textstyle}
\newcommand{\dsps}{\displaystyle}
\newcommand{\fnz}{\footnotesize}
\newcommand{\scz}{\scriptsize}
\newcommand{\be}{
\begin{equation}
}
\newcommand{\bel}[1]{
\begin{equation}
\label{#1}}
\newcommand{\ee}{
\end{equation}
}
\newcommand{\eqnl}[2]{
\begin{equation}
\label{#1}{#2}
\end{equation}
}
\newtheorem{subn}{\name}
\renewcommand{\thesubn}{}
\newcommand{\bsn}[1]{\def\name{#1}
\begin{subn}}
\newcommand{\esn}{
\end{subn}}
\newtheorem{sub}{\name}[section]
\newcommand{\dn}[1]{\def\name{#1}}   
\newcommand{\bs}{
\begin{sub}}
\newcommand{\es}{
\end{sub}}
\newcommand{\bsl}[1]{
\begin{sub}\label{#1}}
\newcommand{\bth}[1]{\def\name{Theorem}
\begin{sub}\label{t:#1}}
\newcommand{\blemma}[1]{\def\name{Lemma}
\begin{sub}\label{l:#1}}
\newcommand{\bcor}[1]{\def\name{Corollary}
\begin{sub}\label{c:#1}}
\newcommand{\bdef}[1]{\def\name{Definition}
\begin{sub}\label{d:#1}}
\newcommand{\bprop}[1]{\def\name{Proposition}
\begin{sub}\label{p:#1}}
\newcommand{\R}{\eqref}
\newcommand{\rth}[1]{Theorem~\ref{t:#1}}
\newcommand{\rlemma}[1]{Lemma~\ref{l:#1}}
\newcommand{\rcor}[1]{Corollary~\ref{c:#1}}
\newcommand{\rdef}[1]{Definition~\ref{d:#1}}
\newcommand{\rprop}[1]{Proposition~\ref{p:#1}}
\newcommand{\BA}{
\begin{array}}
\newcommand{\EA}{
\end{array}}
\newcommand{\BAN}{\renewcommand{\arraystretch}{1.2}
\setlength{\arraycolsep}{2pt}
\begin{array}}
\newcommand{\BAV}[2]{\renewcommand{\arraystretch}{#1}
\setlength{\arraycolsep}{#2}
\begin{array}}
\newcommand{\BSA}{
\begin{subarray}}
\newcommand{\ESA}{
\end{subarray}}
\newcommand{\BAL}{
\begin{aligned}}
\newcommand{\EAL}{
\end{aligned}}
\newcommand{\BALG}{
\begin{alignat}}
\newcommand{\EALG}{
\end{alignat}}
\newcommand{\BALGN}{
\begin{alignat*}}
\newcommand{\EALGN}{
\end{alignat*}}
\newcommand{\note}[1]{\textit{#1.}\hspace{2mm}}
\newcommand{\Proof}{\note{Proof}}
\newcommand{\qeda}{\hspace{10mm}\hfill $\square$}
\newcommand{\qed}{\\
${}$ \hfill $\square$}
\newcommand{\Remark}{\note{Remark}}
\newcommand{\modin}{$\,$\\
[-4mm] \indent}
\newcommand{\forevery}{\quad \forall}
\newcommand{\set}[1]{\{#1\}}
\newcommand{\setdef}[2]{\{\,#1:\,#2\,\}}
\newcommand{\setm}[2]{\{\,#1\mid #2\,\}}
\newcommand{\lra}{\longrightarrow}
\newcommand{\lla}{\longleftarrow}
\newcommand{\llra}{\longleftrightarrow}
\newcommand{\Lra}{\Longrightarrow}
\newcommand{\Lla}{\Longleftarrow}
\newcommand{\Llra}{\Longleftrightarrow}
\newcommand{\warrow}{\rightharpoonup}
\newcommand{
\paran}[1]{\left (#1 \right )}
\newcommand{\sqbr}[1]{\left [#1 \right ]}
\newcommand{\curlybr}[1]{\left \{#1 \right \}}
\newcommand{\abs}[1]{\left |#1\right |}
\newcommand{\norm}[1]{\left \|#1\right \|}
\newcommand{
\paranb}[1]{\big (#1 \big )}
\newcommand{\lsqbrb}[1]{\big [#1 \big ]}
\newcommand{\lcurlybrb}[1]{\big \{#1 \big \}}
\newcommand{\absb}[1]{\big |#1\big |}
\newcommand{\normb}[1]{\big \|#1\big \|}
\newcommand{
\paranB}[1]{\Big (#1 \Big )}
\newcommand{\absB}[1]{\Big |#1\Big |}
\newcommand{\normB}[1]{\Big \|#1\Big \|}

\newcommand{\thkl}{\rule[-.5mm]{.3mm}{3mm}}
\newcommand{\thknorm}[1]{\thkl #1 \thkl\,}
\newcommand{\trinorm}[1]{|\!|\!| #1 |\!|\!|\,}
\newcommand{\bang}[1]{\langle #1 \rangle}
\def\angb<#1>{\langle #1 \rangle}
\newcommand{\vstrut}[1]{\rule{0mm}{#1}}
\newcommand{\rec}[1]{\frac{1}{#1}}
\newcommand{\opname}[1]{\mbox{\rm #1}\,}
\newcommand{\supp}{\opname{supp}}
\newcommand{\dist}{\opname{dist}}
\newcommand{\myfrac}[2]{{\displaystyle \frac{#1}{#2} }}
\newcommand{\myint}[2]{{\displaystyle \int_{#1}^{#2}}}
\newcommand{\mysum}[2]{{\displaystyle \sum_{#1}^{#2}}}
\newcommand {\dint}{{\displaystyle \int\!\!\int}}
\newcommand{\q}{\quad}
\newcommand{\qq}{\qquad}
\newcommand{\hsp}[1]{\hspace{#1mm}}
\newcommand{\vsp}[1]{\vspace{#1mm}}
\newcommand{\ity}{\infty}
\newcommand{\prt}{\partial}
\newcommand{\sms}{\setminus}
\newcommand{\ems}{\emptyset}
\newcommand{\ti}{\times}
\newcommand{\pr}{^\prime}
\newcommand{\ppr}{^{\prime\prime}}
\newcommand{\tl}{\tilde}
\newcommand{\sbs}{\subset}
\newcommand{\sbeq}{\subseteq}
\newcommand{\nind}{\noindent}
\newcommand{\ind}{\indent}
\newcommand{\ovl}{\overline}
\newcommand{\unl}{\underline}
\newcommand{\nin}{\not\in}
\newcommand{\pfrac}[2]{\genfrac{(}{)}{}{}{#1}{#2}}

\def\ga{\alpha}     \def\gb{\beta}       \def\gg{\gamma}
\def\gc{\chi}       \def\gd{\delta}      \def\ge{\epsilon}
\def\gth{\theta}                         \def\vge{\varepsilon}
\def\gf{\phi}       \def\vgf{\varphi}    \def\gh{\eta}
\def\gi{\iota}      \def\gk{\kappa}      \def\gl{\lambda}
\def\gm{\mu}        \def\gn{\nu}         \def\gp{\pi}
\def\vgp{\varpi}    \def\gr{\rho}        \def\vgr{\varrho}
\def\gs{\sigma}     \def\vgs{\varsigma}  \def\gt{\tau}
\def\gu{\upsilon}   \def\gv{\vartheta}   \def\gw{\omega}
\def\gx{\xi}        \def\gy{\psi}        \def\gz{\zeta}
\def\Gg{\Gamma}     \def\Gd{\Delta}      \def\Gf{\Phi}
\def\Gth{\Theta}
\def\Gl{\Lambda}    \def\Gs{\Sigma}      \def\Gp{\Pi}
\def\Gw{\Omega}     \def\Gx{\Xi}         \def\Gy{\Psi}

\def\CS{{\mathcal S}}   \def\CM{{\mathcal M}}   \def\CN{{\mathcal N}}
\def\CR{{\mathcal R}}   \def\CO{{\mathcal O}}   \def\CP{{\mathcal P}}
\def\CA{{\mathcal A}}   \def\CB{{\mathcal B}}   \def\CC{{\mathcal C}}
\def\CD{{\mathcal D}}   \def\CE{{\mathcal E}}   \def\CF{{\mathcal F}}
\def\CG{{\mathcal G}}   \def\CH{{\mathcal H}}   \def\CI{{\mathcal I}}
\def\CJ{{\mathcal J}}   \def\CK{{\mathcal K}}   \def\CL{{\mathcal L}}
\def\CT{{\mathcal T}}   \def\CU{{\mathcal U}}   \def\CV{{\mathcal V}}
\def\CZ{{\mathcal Z}}   \def\CX{{\mathcal X}}   \def\CY{{\mathcal Y}}
\def\CW{{\mathcal W}} \def\CQ{{\mathcal Q}}
\def\BBA {\mathbb A}   \def\BBb {\mathbb B}    \def\BBC {\mathbb C}
\def\BBD {\mathbb D}   \def\BBE {\mathbb E}    \def\BBF {\mathbb F}
\def\BBG {\mathbb G}   \def\BBH {\mathbb H}    \def\BBI {\mathbb I}
\def\BBJ {\mathbb J}   \def\BBK {\mathbb K}    \def\BBL {\mathbb L}
\def\BBM {\mathbb M}   \def\BBN {\mathbb N}    \def\BBO {\mathbb O}
\def\BBP {\mathbb P}   \def\BBR {\mathbb R}    \def\BBS {\mathbb S}
\def\BBT {\mathbb T}   \def\BBU {\mathbb U}    \def\BBV {\mathbb V}
\def\BBW {\mathbb W}   \def\BBX {\mathbb X}    \def\BBY {\mathbb Y}
\def\BBZ {\mathbb Z}

\def\GTA {\mathfrak A}   \def\GTB {\mathfrak B}    \def\GTC {\mathfrak C}
\def\GTD {\mathfrak D}   \def\GTE {\mathfrak E}    \def\GTF {\mathfrak F}
\def\GTG {\mathfrak G}   \def\GTH {\mathfrak H}    \def\GTI {\mathfrak I}
\def\GTJ {\mathfrak J}   \def\GTK {\mathfrak K}    \def\GTL {\mathfrak L}
\def\GTM {\mathfrak M}   \def\GTN {\mathfrak N}    \def\GTO {\mathfrak O}
\def\GTP {\mathfrak P}   \def\GTR {\mathfrak R}    \def\GTS {\mathfrak S}
\def\GTT {\mathfrak T}   \def\GTU {\mathfrak U}    \def\GTV {\mathfrak V}
\def\GTW {\mathfrak W}   \def\GTX {\mathfrak X}    \def\GTY {\mathfrak Y}
\def\GTZ {\mathfrak Z}   \def\GTQ {\mathfrak Q}

\font\Sym= msam10 
\def\SYM#1{\hbox{\Sym #1}}
\newcommand{\bdw}{\prt\Gw\xspace}
\medskip
\mysection {Introduction}

Let  $\Gw$ be a bounded domain in $\mathbb R^{N}$ with smooth boundary and $\gm$ a Radon measure on $\prt\Gw$. In this paper we consider first the problem of finding a function $u$ solution of
\begin {equation}
\label {exp1} -\Gd u+e^{u}-1=0
\end {equation}
in $\Gw$ satisfying $u=\gm$ on $\prt\Gw$. Let $\gr (x)=\dist (x,\prt \Gw)$, then this problem admits a {\bf weak formulation}: {\it find a function $u\in L^1(\Gw)$ such that 
$e^{u}\in L^1_\gr(\Gw)$ satisfying }
\begin {equation}
\label {exp1w} \myint{\Gw}{}\left(-u\Gd\gz+(e^{u}-1)\gz\right)dx=-\myint{\prt\Gw}{}\frac{\prt\gz}{\prt\gn}d\gm\qquad\forall\gz\in W^{1,\infty}_0(\Gw)\cap W^{2,\infty}(\Gw). 
\end {equation}
This equation has been initiated by Grillot and V\'eron \cite {GrV} in 2-dim in the framework of the boundary trace theory.  Much works on boundary trace problems for equation of the type 
 \begin {equation}
\label {pow1} -\Gd u+u^q=0
\end {equation}
with $q>1$),  have been developed by Le Gall \cite {LG}, Marcus and V\'eron \cite{MV1}, \cite{MV2}, Dynkin and Kuznetsov \cite{DK1}, \cite{DK2},  respectively by purely probabilistic methods, by purely analytic methods or by a combination of the preceding aspects. One of the main features of the problem with power nonlinearities is the existence of a critical exponent $q_c=\frac{N+1}{N-1}$ which is linked to the existence of boundary removable sets. Existence of boundary removable points have been discovered by Gmira and V\'eron \cite {GmV}. Let us recall briefly the main results for ($\ref{pow1}$):\smallskip

\noindent (i) If $1<q<q_c$, then for any $\gm\in\mathfrak M_+(\prt\Gw)$ there exists a unique function $u\in L^1_+(\Gw)\cap L^q_\gr(\Gw)$ which satisfies ($\ref{pow1}$) in $\Gw$ and takes the value $\gm$ on $\prt\Gw$ in the following weak sense
\begin {equation}
\label {poww} \myint{\Gw}{}\left(-u\Gd\gz+u^q\gz\right)dx=-\myint{\prt\Gw}{}\frac{\prt\gz}{\prt\gn}d\gm\qquad\forall\gz\in W^{1,\infty}_0(\Gw)\cap W^{2,\infty}(\Gw). 
\end {equation}

\noindent (ii) If $q\geq q_c$, the above problem can be solved if and only if $\gm$ vanishes on boundary Borel subsets with zero $C_{\frac{2}{q},q'}$-Bessel capacity. Furthermore a boundary compact set is removable if and only if it has zero $C_{\frac{2}{q},q'}$-capacity.\\

In this article we adapt some of the ideas used for $(\ref{pow1})$ to problem
\begin {equation}
\label {exp2}\BA {ll} -\Gd u+e^{u}-1=0\qquad&\mbox {in }\;\Gw\\
\phantom{-\Gd +e^{u}-1}
u=\gm\qquad&\mbox {on }\;\prt\Gw.
\EA
\end {equation}
Following the terminology of \cite{BMP} we say that a measure $\gm\in \mathfrak M(\prt\Gw)$ is {\bf good} if $(\ref{exp2})$ admits a weak solution. Let  $P^\Gw(x,y)$ (resp. $G^\Gw(x,y)$) be the Poisson kernel (resp. the Green kernel) in $\Gw$ and $\mathbb P^\Gw[\gm]$ 
the Poisson potential of a  boundary mesure $\gm$ (resp. 
 $\mathbb G^\Gw[\phi]$ the Green potential of a bounded measure $\phi$ defined in $\Gw$). A boundary measure $\gm$ which satisfies
 \begin {equation}
\label {exp4}
\quad\exp(\mathbb P^\Gw[\gm])\in L^{1}(\Gw;\gr dx).
\end {equation}
 is called {\bf admissible}. Since $\mathbb P^\Gw[\gm]$ is a supersolution for $(\ref{exp1})$, an admissible measure is good (see \cite{Ve}). 
 Our first result which extends a previous one obtained in \cite{GrV} is the following.\\

\noindent{\bf Theorem A}. {\it Suppose $\gm\in \mathfrak M(\prt\Gw)$ admits Lebesgue decomposition $\gm=\gm_S+\gm_R$ where $\gm_S$ and $\gm_R$ are mutually singular and $\gm_R$ is absolutely continuous with respect to the (N-1)-dim Hausdorff measure $dH^{N-1}$. If
\begin {equation}
\label {exp3}
\quad\exp(\mathbb P^\Gw[\gm_{S}])\in L^{1}(\Gw;\gr dx),
\end {equation}
then $\gm$ is good.}\\

In order to go further in the study of good measures, it is necessary to introduce an {\it Orlicz capacity} modelized on the Legendre transform of  $r\mapsto p(r):=e^{r}-1$. These capacities have been studied by Aissaoui and Benkirane \cite {AsBe} and they inherit  most of the properties of the Bessel capacities. The capacity $C_{N^{L\ln L}}$ associated to the problem is constructed later and it has strong connexion with Hardy-Littlewood maximal function. In this capacity framework  we obtain the following types of results:\\

\noindent{\bf Theorem B}. {\it Let $\gm\in\mathfrak M_+(\prt\Gw)$ be a good measure, then $\gm$ vanishes on boundary Borel subsets $E$ with zero $C_{N^{L\ln L}}$-capacity.}\\

We also give below a result of removability of boundary singularities.\\

\noindent{\bf Theorem C}. {\it Let $K\subset\prt\Gw$ be a compact subset with zero $C_{N^{L\ln L}}$-capacity. Suppose $u\in C(\overline\Gw\setminus K)\cap C^2(\Gw)$ is a positive solution of $(\ref{exp1})$ in $\Gw$ which vanishes on $K$, then $u$ is identically zero.}\\

In the last part of this paper we apply  this approach to the problem
\begin {equation}
\label {exp11} -\Gd u+e^{u}-1=\gm,
\end {equation}
where $\gm$ is a bounded measure, as well as removability questions for internal singularities of solutions of $(\ref{exp1})$. In that case the natural capacity associated to the problem is 
\begin {equation}\label {exp12} 
C_{\Gd^{L\ln L}}(K)=\inf\{\norm{M[\Gd\eta]}_{L^1}:\eta\in C^2_0(\Gw):0\leq \eta\leq 1,\eta=1\text{ in a neighborhood of} K\}
\end {equation}
where $M[.]$ denotes Hardy-Littlewood's maximal function.\medskip

\noindent{\bf Theorem D}. {\it Let $\gm\in\mathfrak M^b_+(\Gw)$ be a bounded good measure, then $\gm$ vanishes on boundary Borel subsets $E$ with zero $C_{\Gd^{L\ln L}}$-capacity.}\medskip

A charaterization of positive measures which have the property of vanishing on Borel subsets $E$ with zero $C_{N^{L\ln L}}$-capacity is also provided. We also give below a result of removability of boundary singularities.\\

\noindent{\bf Theorem E}. {\it Let $K\subset\Gw$ be a compact subset with zero $C_{\Gd^{L\ln L}}$-capacity. Suppose $u\in C(\Gw\setminus K)\cap C^2(\Gw)$ is a positive solution of $(\ref{exp1})$ in $\Gw\setminus K$ which vanishes on $\prt\Gw$, then $u$ is identically zero.}\medskip\\

\nind { This note is essentially derived from the preliminary report \cite {Ve1}, written in 2004 and left escheated since this period. }


\mysection {Good measures}

\noindent{\it Proof of Theorem A.}  For $k>0$, 
set $\gm_{R,k}=\inf\{k,\gm_{R}\}$ 
and denote by $u_{k}$ the 
solution of 
\begin {equation}
\label {expk}\BA {ll} -\Gd u_{k}+e^{u_{k}}-1=0\qquad&\mbox {in }\;\Gw\\\phantom{-\Gd +e^{u_{k}}-1}
 u_{k}=\gm_{S}+\gm_{R,k}\qquad&\mbox {on }\;\prt\Gw.
\EA
\end {equation}
Such a solution exists because 
$$
\exp(\mathbb P^\Gw[\gm_{S}+\gm_{R,k}])\leq e^k
\exp(\mathbb P^\Gw[\gm_{S}]) $$
by the maximum principle, and ($\ref {exp3}$) implies that 
$\exp(\mathbb P^\Gw[\gm_{S}+\gm_{R,k}])-1\in L^{1}(\Gw;\gr dx)$. The 
sequence $u_{k}$ is nondecreasing. Since, for any $\gz\in 
C_{c}^{1,1}(\bar\Gw)$, 
$$
\int_{\Gw}(-u_{k}\Gd\gz+(e^{u_{k}}-1)\gz)dx=\int_{\prt\Gw}
\myfrac {\prt\gz}{\prt\gn}d(\gm_{S}+\gm_{R,k}), $$
if we take in particular for test function $\gz$ the solution $\gz_{0}$ of 
\begin {equation}
\label {expk}\BA {ll} -\Gd \gz_{0}=1\qquad&\mbox {in }\;\Gw\\\phantom{-\Gd }
 \gz_{0}=0\quad&\mbox {on }\;\prt\Gw,
\EA
\end {equation}
we get
\begin {equation}
\label {expk'}
\int_{\Gw}(u_{k}+(e^{u_{k}}-1)\gz_{0})dx=-\int_{\prt\Gw}
\myfrac {\prt\gz_{0}}{\prt\gn}d(\gm_{S}+\gm_{R,k})\leq c\norm 
{\gm}_{\mathfrak M}.
\end {equation}
Thus $u=\lim_{k\to\ity}u_{k}$ is integrable, $$
\int_{\Gw}(u+(e^{u}-1)\gz_{0})dx\leq c\norm 
{\gm}_{\mathfrak M}, $$
and the convergence of $u_{k}$ and $e^{u_{k}}$ to $u$ and $e^{u}$ hold respectively in
$L^{1}(\Gw)$ and $L^{1}(\Gw;\gr dx)$ and $u$ satisfies $(\ref {exp1w})$.\qeda \\

The proof of the next result is directly inspired by \cite {BMP} where nonlinear Poisson equations are treated.\\

\bprop{2} The following properties hold:\smallskip

\nind (i) If $\gm\in \mathfrak M_{+}(\prt\Gw)$ is a  good measure, then
any $\tilde\gm\in \mathfrak M_+(\prt\Gw)$ smaller than $\gm$ is good.\smallskip

\nind (ii) Let $\{\gm_{n}\}$ be an increasing sequence of good measures
 which converges to $\gm$ in the weak sense of measures. Then 
 $\gm$ is good.
 \smallskip

\nind (iii)  If $\gm\in \mathfrak M_{+}(\prt\Gw)$ is a  good measure and $f\in L_{+}^{1}(\prt\Gw)$, then $f+\gm$ is a good measure. \es

\nind\Proof (i) Let $u=u_{\gm}$ be the solution of $(\ref{exp2})$ and 
$w=\inf\{u,\mathbb P^{\Gw}[\tilde\gm]\}$. Since $\mathbb 
P^{\Gw}[\tilde\gm]$ is a supersolution for $(\ref{exp1})$, $w$ is a 
supersolution too. Furthermore $w$ is nonnegative and $e^{w}-1\in 
L^{1}(\Gw;\gr dx)$. By Doob's theorem $w$ admits a boundary trace 
$\gm^{*}\in\mathfrak M_{+}(\prt\Gw)$ and $\gm^{*}\leq \tilde\gm\leq\gm$. Let $w^{*}$ be the
solution of 
$$
\BA {ll} -\Gd w^{*}+e^u-1=0\qquad&\mbox {in }\;\Gw\\\phantom{-\Gd +e^u-1}
 w^{*}=\tilde\gm\quad&\mbox {on }\;\prt\Gw.
\EA$$
then $u\geq w\geq w^{*}$ and \cite {MV3},  
$$
\lim_{t}\int_{\prt\Gw_{t}}w^{*}(t,.)\eta dS_{t}=\int_{\prt\Gw}\eta 
d\tilde\gm\forevery\eta\in C(\prt\Gw), $$
(here we denote by $\prt\Gw_{t}$ the set of $x\in\Gw$ such that 
$\rho (x)=t>0$). This implies that the boundary trace of $w^{*}$ is 
$\tilde\gm$ and thus $\gm^{*}=\tilde\gm$. Set $\Gw_{t}=\{x\in\Gw: \rho 
(x)>t\}$ and let $v_{t}$ we the solution of 
$$
\BA {ll} -\Gd v_{t}+e^{v_{t}}-1=0\qquad&\mbox {in }\;\Gw_{t}\\\phantom{-\Gd +e^{v_{t}}-1}
v_{t}=w\quad&\mbox {on }\;\prt\Gw_{t}.
\EA$$
Then $v_{t}\leq w$ in $\Gw_{t}$. Furthermore $$
0<t'<t\Longrightarrow v_{t'} \leq v_{t}\quad\mbox {in }\Gw_{t}.$$
Then $\tilde u=\lim_{t\to 0}v_{t}$ exists, the convergence holds in 
$L^{1}(\Gw)$ and $e^{v_{t}}\to e^{\tilde u}$ in $L^{1}(\Gw;\gr dx)$ 
(here we use the fact that $e^w\in L^{1}(\Gw;\gr dx)$. Because 
$$
\lim_{t\to 0}\int_{\prt\Gw_{t}}\tilde w(t,.)\eta dS_{t}=\int_{\prt\Gw}\eta 
d\tilde\gm\forevery\eta\in C(\prt\Gw), $$
and $v_{t}=\tilde w$ on $\prt\Gw_{t}$, is follows that $\tilde u$ admits 
$\tilde\gm$ for boundary trace and thus $\tilde 
u=u_{\tilde\gm}$.\smallskip 

\nind (ii)  Let $u_{n}=u_{\gm_{n}}$ be the solutions of $(\ref {exp2})$ 
with boundary value $\gm_{n}$. The sequence $\{u_{n}\}$ is 
increasing. Since 
\begin {equation}
\label {expk'}
\int_{\Gw}(u_{n}+(e^{u_{n}}-1)\gz_{0})dx=-\int_{\prt\Gw}
\myfrac {\prt\gz_{0}}{\prt\gn}d\gm_{n}\leq -\int_{\prt\Gw}
\myfrac {\prt\gz_{0}}{\prt\gn}d\gm,
\end {equation}
we 
conclude as in the proof of Theorem 1, 
that $u_{n}$ increases and converges to a solution $u=u_{\gm}$ of 
$(\ref {exp2})$ with  boundary value $\gm$.\smallskip

\nind (iii) In the proof of (i) we have actually used the following 
result : {\it Let $w$ be a nonnegative supersolution of $(\ref {exp1})$
such that $e^w\in L^{1}(\Gw;\gr dx)$ and let $\gm\in \mathfrak 
M_{+}(\prt\Gw)$ be the boundary trace of $w$. Then $\gm$ is good}. Let $f\in
L^{1}_{+}(\prt\Gw)$ and $\gm$ be an good measure. We 
denote by $u=u_{\gm}$ the solution of $(\ref {exp2})$. For $k>0$, set $f_{k}=\min\{k,f\}$. 
The function $w_{k}=u_{\gm}+\mathbb P^{\Gw}[f_{k}]$ is a nonnegative 
supersolution, and, since $\mathbb P^{\Gw}[f_{k}]\leq k$, 
$e^{w_{k}}\in L^{1}(\Gw;\gr dx)$. Furthermore the boundary trace of 
$w_{k}$ is $\gm+f_{k}$. Therefore $\gm+f_{k}$ is good. We 
conclude by II that $\gm+f$ is good\qeda\medskip

\nind \Remark The assertions (i) and (ii) in Theorem 1 are still valid if 
we replace $r\mapsto e^r-1$ by any continuous nondecreasing function 
$f$ vanishing at $0$.

\mysection {The Orlicz space framework} 
\subsection{Orlicz capacities}

The set of nonnegative  measures $\gm$ on $\prt\Gw$
such that 
\begin {equation}
\label {orl1}
\exp (\mathbb P^{\Gw}[\gm])\in L^{1}(\Gw;\gr dx)
\end {equation}
 is not a linear space, but it is a convex 
subset of $\mathfrak M_{+}(\prt\Gw)$. The role of this set comes from the fact that any measure in $M^{exp}(\prt\Gw)$ is good. Put 
$$
p(t)={\rm sgn} (s) (e^s-1),\;P(t)=e^{\abs t}-1-{\abs t},$$
and $$
\bar p(s)={\rm sgn} (s) \ln ({\abs s}+1),
\; P^{*}(t)=({\abs t}+1)\ln ({\abs t}+1)-{\abs t}.$$
Then $P$ and $P^{*}$ are complementary functions in the sense of Legendre. Furthermore 
Young inequality holds 
$$
xy\leq P(x)+P^{*}(y)\quad\forevery (x,y)\in \mathbb R\ti \mathbb R, $$
with equality if and only if $x=\bar p(y)$ or $y=p(x)$. It is 
classical to define
\begin {equation}
\label {orlP} M_{P}(\Gw;\gr dx)=\{\phi\in L^{1}_{loc}(\Gw): P(\phi)\in L^{1}(\Gw;\gr 
dx)\},
\end {equation}
\begin {equation}
\label {orlP^{*}} M_{P^{*}}(\Gw;\gr dx)=\{\phi\in L^{1}_{loc}(\Gw): P^{*}(\phi)\in
L^{1}(\Gw;\gr dx)\}.
\end {equation}
The Orlicz spaces $L_{P}(\Gw;\gr dx)$ and $L_{P^{*}}(\Gw;\gr dx)$ are 
the vector spaces spanned respectively by $M_{P}(\Gw;\gr dx)$ and $M_{P^{*}}(\Gw;\gr 
dx)$. They are endowed with the Luxemburg norms
\begin {equation}
\label {orlP1}
\norm \phi_{L_{P_\gr}}=
\inf\left\{k>0:\int_{\Gw}P\left(\myfrac {\gf}{k}\right)\gr dx\leq 1\right\}.
\end {equation}
and
\begin {equation}
\label {orlP^{*}1}
\norm \phi_{L_{P^{*}_\gr}}=
\inf\left\{k>0:\int_{\Gw}P^{*}\left(\myfrac {\gf}{k}\right)\gr dx\leq 1\right\}.
\end {equation}
Furthermore the H\" older-Young inequality asserts
\cite {KrRu}
\begin {equation}
\label {YH}
\abs{\int_{\Gw}\phi\,\psi\,\gr\,dx}\leq \norm \phi_{L_{P_\gr}}\norm \psi_{L_{P^{*}_\gr}}
\forevery (\phi,\psi)\in L_{P}(\Gw;\gr dx)\ti L_{P^{*}}(\Gw;\gr dx).
\end {equation}
Since $P^{*}$ satisfies the $\Gd_{2}$-condition, $M_{P^{*}}(\Gw;\gr 
dx)=L_{P^{*}}(\Gw;\gr dx)$ and $L_{P}(\Gw;\gr dx)$ is the dual space of $L_{P^{*}}(\Gw;\gr dx)$, 
(see  \cite {FuSe}, \cite{AsBe}). Furthermore, since $$
\myfrac {\abs a\ln (1+\abs a)}{2}\leq P^{*}(a)\leq \abs a\ln (1+\abs a)
\forevery a\in \mathbb R, $$
 the space $L_{P^{*}}(\Gw;\gr dx)$ 
is associated with the class $L\ln L(\Gw;\gr dx)$ and to the Hardy-Littlewood maximal
function (see \cite{FuSe}). We recall its definition: we 
consider a cube $Q_{0}$ containing $\bar\Gw$, with sides parallel to the 
axes. If $f\in L^{1}(\Gw)$ we denote by $\tilde f$ its extension by 
$0$ in $Q_{0}\setminus \Gw$ and put 
$$
M_{Q_{0}}[f](x)=\sup\left\{\myfrac{1}{\abs Q}\int_{Q}\abs f(y)dy:Q\in 
{\cal Q}_{x}\right\} $$
where ${\cal Q}_{x}$ denotes the set of all cubes containing $x$ and contained in $Q_{0}$,
with sides parallel to the axes. Thus 
\begin {equation}\label{HL}
\norm f_{L\ln L_{\gr}}:= \int_{Q_{0}}M_{Q_{0}}[f](x) \gr dx\approx\norm f_{L_{P^{*}_\gr}}.
\end{equation}

\bdef{bdrym} The space of all measures on $\prt\Gw$ such that $\mathbb P^\Gw[\gm]\in L_{P}(\Gw;\gr dx) $
is denoted by $B^{\exp}(\prt\Gw)$ and endowed with the norm 
\begin {equation}
\label {Bexp}
\norm {\gm}_{B^{\exp}}=\norm {\mathbb P^\Gw[\gm]}_{L_{P_{\gr}}}.
\end{equation}
The subset of measures such that $\exp(\mathbb P^\Gw[\gm])\in L^{1}(\Gw;\gr dx) $
is denoted by $\mathfrak M^{\exp}(\prt\Gw)$.
\es

The following result follows immediately from the definition of the Luxemburg norm.

\bprop{0} If $\gm\in B^{\exp}(\prt\Gw)$ there exists $a_0>0$ such that $a\gm\in\mathfrak M^{\exp}(\prt\Gw)$ for all $0\leq a<a_0$. Conversely, if $\gm\in\mathfrak M^{\exp}(\prt\Gw)$, then $a\gm\in B^{\exp}(\prt\Gw)$ for all $a>0$.
\es

The analytic charaterization of $B^{\exp}(\prt\Gw)$ can be done in 
introducing the space of normal derivatives of Green potentials of $L\ln L$ functions:
\begin {equation}
\label {NLlnL} N^{L\ln L}(\prt\Gw)=
\left\{\eta: \gr^{-1}\Gd(\gr^{*}\mathbb P^\Gw[\eta])\in L\ln 
L(\Gw;\gr dx)\right\}.
\end{equation}
where $\gr^{*}$ is a smooth positive function with value $\gr$ in a 
neighborhood of $\prt\Gw$. Then
\begin {equation}
\label {dual}
\abs{\int_{\prt\Gw}\eta d\gm}=\abs{\int_{\Gw}\mathbb 
P^\Gw[\gm]\Gd(\gr^{*}\mathbb P^\Gw[\eta])\,dx}
\leq \norm{\mathbb P^\Gw[\gm]}_{L_{P_{\gr}}}
\norm{\gr^{-1}\Gd(\gr^{*}\mathbb P^\Gw[\eta])}_{L_{P^{*}_{\gr}}}.
\end{equation}

We take for norm on $N^{L\ln L}(\prt\Gw)$
\begin {equation}
\label {NLlnL2}
\norm {\eta}_{N^{L\ln L}}=
\norm{\gr^{-1}\Gd(\gr^{*}\mathbb P^\Gw[\eta])}_{L_{P^{*}_{\gr}}},
\end{equation}
and define the $C_{N^{L\ln L}}$-capacity of a compact subset $K$ of 
$\prt\Gw$ by
\begin {equation}
\label {NLlnL3} C_{N^{L\ln L}}(K)=\inf\{\norm {\eta}_{N^{L\ln L}}:\eta\in C^{2}(\prt\Gw), 
0\leq\eta\leq 1,\eta\geq 1\mbox { in a neighborhood of }K\}.
\end{equation}
Considering the bilinear form $\CH$ on $L_{P^{*}_{\gr}}(\prt\Gw)\ti L_{P_{\gr}}(\prt\Gw)$ 
\begin {equation}
\label {NLlnL4} 
\CH (\eta,\gm):=-\int_{\Gw}\mathbb 
P^\Gw[\gm]\Gd(\gr^{*}\mathbb P^\Gw[\eta])\,dx
\end{equation}
then
\begin {equation}
\label {NLlnL5} \BA {l}
\CH (\eta,\gm)=-\myint{\Gw}{}\myint{\prt\Gw}{}P^\Gw(x,y)d\gm(y)\Gd(\gr^{*}\mathbb P^\Gw[\eta])(x)\,dx\\[4mm]
\phantom{\CH (\eta,\gm)}
=-\myint{\prt\Gw}{}\myint{\Gw}{}\Gd(\gr^{*}\mathbb P^\Gw[\eta])(x)P^\Gw(x,y)\,dx\,d\gm(y).
\EA\end{equation}
It is classical to define 
\begin {equation}
\label {NLlnL6} \BA {l}C^*_{N^{L\ln L}}(K)=\sup\{\gm (K):\gm\in \mathfrak M_+(\prt\Gw), 
\gm(K^c)=0,\norm{\BBP^{\Gw}[\gm]}_{L_{P_\gr}}\leq 1\}.
\EA\end{equation}
The following result due to Fuglede \cite{Fug} (and to Aissaoui-Benkirane in the Orlicz space framework \cite {AsBe}) is a consequence of the Kneser-Fan min-max theorem.
\bprop{equ} For any compact set $K\subset\prt\Gw$, there holds
\begin {equation}\label {NLlnL7} 
C^*_{N^{L\ln L}}(K)=C_{N^{L\ln L}}(K).
\end{equation}
\es

As a direct consequence of $(\ref {dual})$, we have the following

\bprop{1} If $\gm\in B_{+}^{\exp}(\prt\Gw) $, it does 
not charge Borel subsets with $C_{N^{L\ln L}}$-capacity zero. \es


\subsection{Good measures and removable sets}


\nind{\it Proof of Theorem B. } Let $K$ be a compact subset with $C_{N^{L\ln L}}$-capacity zero. There exist a
sequence $\{\eta_n\}\subset C^2(\prt\Gw)$ such that $0\leq \eta_n \leq 1$, $\eta_n=1$ in a
neighborhood of $K$ and 
\begin {equation}
\label {AM1}
\lim_{n\to\infty}\norm {\eta_n}_{N^{L\ln L}}=
\norm {\gr^{-1}\Gd(\gr^*\mathbb P^\Gw[\eta_n])}_{L_{P^*_\gr}}=0.
\end {equation}
Take $\gr^*\mathbb P^\Gw[\eta_n])$ as a test function, then 
$$
\int_\Gw\left(-u\Gd(\gr^*\mathbb P^\Gw[\eta_n])+(e^{u}-1)\gr^*\mathbb P^\Gw[\eta_n]))
\right)dx =-\int_{\prt\Gw} \myfrac {\prt(\gr^*\mathbb P^\Gw[\eta_n]))}{\prt \gn}d\gm $$
Since $\myfrac {\prt(\gr^*\mathbb P^\Gw[\eta_n]))}{\prt \gn}=\eta_n $ and $\gm>0$, there
holds $-\myint{\prt\Gw}{} \myfrac {\prt(\gr^*\mathbb P^\Gw[\eta_n]))}{\prt \gn}d\gm\geq \gm
(K). $ Furthermore
\begin {equation}
\label {AM2}
\abs {\int_\Gw u\Gd(\gr^*\mathbb P^\Gw[\eta_n]) dx}\leq 
\norm {u}_{L_{P_\gr}}\norm {\gr^{-1}\Gd(\gr^*\mathbb P^\Gw[\eta_n])}_{L_{P^*_\gr}}.
\end {equation}
Then $$
\gm (E)\leq \int_\Gw (e^{u}-1)\gr^*\mathbb P^\Gw[\eta_n]) dx+
\norm {u}_{L_{P_\gr}}\norm {\gr^{-1}\Gd(\gr^*\mathbb P^\Gw[\eta_n])}_{L_{P^*_\gr}}. $$
By the same argument as in \cite{BMP}, 
$\lim_{n\to\ity} \gr^*\mathbb P^\Gw[\eta_n]=0$, a.e. in $\Gw$, and there exists a
nonnegative 
$L^1_\gr$-function $\Gf$ such that $0\leq \gr^*\mathbb P^\Gw[\eta_n]\leq \Gf$. By $(\ref
{AM1})$, $(\ref {AM2})$ and 
Lebesgue's theorem, $\gm (E)=0$.\qeda 
\bdef {remov}A subset $E\subset\prt\Gw$ is said removable for equation $(\ref{exp1})$, if and only if any positive
solution $u\in C^2(\Gw)$ of $(\ref {exp1})$ in $\Gw$, which is continuous in $\overline\Gw\setminus E$ and
vanishes on $\prt\Gw\setminus E$, is identically zero.
\es\medskip


\nind{\it Proof of Theorem C.} Let $u\in C(\Bar\Gw\setminus K)$ be a solution of $(\ref {exp1})$ which is zero
on $\prt\Gw\setminus K$. Let $\{\eta_n\}\subset C^2(\prt\Gw)$ such that $0\leq \eta_n\leq1$,
$\eta_n=1$ in a neighborhood $\CV$ of $K$ and $(\ref {AM1})$ holds. Put $\gth_n =1-\-
\eta_n$. Put $\gr_K(x)=\dist (x,K)$. Then, as a consequence of Keller-Osserman estimate
and the fact that $u$ vanishes on $K^c$, there holds $$
u(x)\leq C\myfrac {\gr(x)\ln (2/\gr_K (x))}{\gr_K (x)}+D. $$
Thus the function $\gz_n=\gr^*\mathbb P^\Gw[\theta_n]$ is an admissible test function for
$u$, and $$
\int_\Gw\left(-u\Gd\gz_n+(e^u-1)\gz_n\right) dx=0. $$
Clearly $\mathbb P^\Gw[\theta_n]=1-\mathbb P^\Gw[\eta_n]$ and 
$$
\Gd\gz_n=\Gd\gr^*-\Gd(\gr^*\mathbb P^\Gw[\eta_n]) $$
Inasmuch we can modify $\gr^*$ in order to have $-\Gd\gr^*\geq 0$, in which case
$\gr^*=\gr$ near $\prt\Gw$ is replaced by $\gr^*\approx \gr$, we derive
\begin {eqnarray*}
-\int_\Gw u\Gd\gz_n\,dx=-\int_\Gw \gz_n^{-1}\Gd\gz_n\,u\gz_ndx
\phantom {-\int_\Gw u\Gd\gz_n\,dx-------------------,,}\\
\geq -2^{- 1}\int_\Gw(e^u-1-u)\gz_n\,dx-\int_\Gw Q(\gz_n^{-1}\Gd(\gr^*\mathbb
P^\Gw[\eta_n]))\,\gz_ndx,
\phantom {---------}
\end {eqnarray*}
where 
$$
Q(r)=(\abs r+2^{-1})\ln (2\abs r+1)-\abs r\leq C\abs r\ln (\abs r+1)\forevery r\in\mathbb
R. $$
Therefore
\begin {eqnarray}
\label {est}
\int_\Gw(e^u-1-u)\gz_n\,dx
\leq 2C\int_\Gw \abs{\Gd(\gr^*\mathbb P^\Gw[\eta_n])}\ln (1+\gr^{-2}\abs{\Gd(\gr^*\mathbb
P^\Gw[\eta_n])})dx,
\end {eqnarray}
since $\gz_n^{-1}\abs{\Gd(\gr^*\mathbb P^\Gw[\eta_n])}\leq \gr^{-2}\abs{\Gd(\gr^*\mathbb
P^\Gw[\eta_n])}$. Furthermore
\begin {eqnarray*}
\ln (1+\gr^{-2}\abs{\Gd(\gr^*\mathbb P^\Gw[\eta_n])})) =-\ln\gr+\ln
(\gr+\gr^{-1}\abs{\Gd(\gr^*\mathbb P^\Gw[\eta_n])}) \\
\leq -\ln\gr+\ln (1+\gr^{-1}\abs{\Gd(\gr^*\mathbb P^\Gw[\eta_n])}) 
\end {eqnarray*}
But (we can assume $\gr\leq 1$)
\begin {eqnarray*}
\int_\Gw \abs{\Gd(\gr^*\mathbb P^\Gw[\eta_n])}\ln (1+\gr^{-2}\abs{\Gd(\gr^*\mathbb
P^\Gw[\eta_n])})dx
\phantom {------------------}\\
\leq 
-\int_\Gw \abs{\Gd(\gr^*\mathbb P^\Gw[\eta_n])}\ln\gr dx +\int_\Gw\abs{\Gd(\gr^*\mathbb
P^\Gw[\eta_n])}
\ln (1+\gr^{-1}\abs{\Gd(\gr^*\mathbb P^\Gw[\eta_n])})dx, 
\end {eqnarray*}
and
\begin {eqnarray*}
\int_\Gw \abs{\Gd(\gr^*\mathbb P^\Gw[\eta_n])}\ln\gr^{-1} dx
\phantom {------------------------------}\\
= \int_{ \{\abs{\Gd(\gr^*\mathbb P^\Gw[\eta_n])}\leq 1\}} \abs{\Gd(\gr^*\mathbb
P^\Gw[\eta_n])}\ln\gr^{-1} dx+
 \int_{\{ \abs{\Gd(\gr^*\mathbb P^\Gw[\eta_n])}> 1\}} \abs{\Gd(\gr^*\mathbb
P^\Gw[\eta_n])}\ln\gr^{-1} dx
 \phantom {---}\\
 \leq \int_{ \{\abs{\Gd(\gr^*\mathbb P^\Gw[\eta_n])}\leq 1\}} \abs{\Gd(\gr^*\mathbb
P^\Gw[\eta_n])}\ln\gr^{-1} dx
 +\int_\Gw\abs{\Gd(\gr^*\mathbb P^\Gw[\eta_n])}
\ln (1+\gr^{-1}\abs{\Gd(\gr^*\mathbb P^\Gw[\eta_n])})dx
\end {eqnarray*}
But $$
\lim_{n\to\infty}\abs{\Gd(\gr^*\mathbb P^\Gw[\eta_n])}=0\quad \mbox {a. e. in }\Gw, $$
at least up to some subsequence. Thus
\begin {eqnarray}
\label {lim1}\lim_{n\to\infty}
\int_\Gw \abs{\Gd(\gr^*\mathbb P^\Gw[\eta_n])}\ln (1+\gr^{-2}\abs{\Gd(\gr^*\mathbb
P^\Gw[\eta_n])})dx=0
\end {eqnarray}
Using $(\ref {est})$, we derive $u=0$.\smallskip

\nind Conversely, assume that $C_{N^{L\ln L}}(K)>0$. By \rprop{equ} there  exists a 
non negative non-zero measure $\gm\in\mathfrak M_+(\prt\Gw)$ such that $\gm (K^c)=0$ in the
space $B_{+}^{exp}(\prt\Gw)$. This means that $\gth\gm\in M_{+}^{exp}(\prt\Gw)$ for some
$\gth>0$. Thus
problem $(\ref {exp2})$ admits a solution.\qeda
\\


{\nind Several questions can be adressed}\medskip

\nind 1- If a measure $\gm$ is good does there exist an
increasing sequence of measures $\{\gm_n\}$ which converges to $\gm$ such that $\gth_n\gm_n$ is admissible 
for some $\gth_n>0$?\medskip

\nind 2-  If a measure $\gm$, singular with respect to $\CH^{N-1}$ is good does it exist an increasing sequence of
admissible measures $\{\gm_n\}$ converging to $\gm$ ?\medskip

\nind 3-  If a measure $\gm$ does not charge Borel sets with $C^{L\ln L}$-capacity zero,
doest it exist $\gth>0$ such that $\gth\gm$ is admissible?\medskip

\nind 4- If a singular measure $\gm$ is good, then $(1-\gd)\gm$ is admissible for any
$\gd\in (0,1)$ ?

\subsection {More general nonlinearities}
\begin {equation}
\label {exp2}\BA {ll} -\Gd u+P(u)=0\qquad&\mbox {in }\;\Gw\\
\phantom{-\Gd +P(u)}
u=\gm\quad&\mbox {on }\;\prt\Gw,
\EA
\end {equation}
where $P$ is a convex increasing function vanishing at $0$ and such that
$\lim_{r\to\infty}P(r)/r=\infty$: In Theorem A-$P$, $(\ref {exp3})$ should be replaced by 
\begin {equation}
\label {f1} P(\mathbb P^\Gw[\gm_{S}])\in L^1(\Gw;\gr\,dx).
\end {equation}
 In \rprop{2}-$P$, (i), (ii) and (iii) still hold. For simplicity we assume that $P$ is a
$N$-function in the 
sense of Orlicz spaces i.e.
 $$
P(r)=\int_{0}^rp(s)ds $$
where $p$ is increasing, vanishes at $0$ and tends to infinity at infinity. Let $P^{*}$ be the 
conjugate $N$-function, $L_{P}(\Gw;\gr\,dx)$ and 
$L_{P^{*}}(\Gw;\gr\,dx)$ the corresponding Orlicz spaces endowed with 
the Luxenburg norms. Then \rprop{1}-$P$ is valid, provided the 
space $$B^{P}(\prt\Gw):=\{\gm\in\mathfrak M(\prt\Gw):\BBP^\Gw[\gm]\in L_{P}(\Gw;\gr\,dx)\}$$ 
endowed with its natural norm replaces $B^{\exp}(\prt\Gw)$ with the norm $(\ref{Bexp})$.
We set 
$$
N^{P^{*}}(\prt\Gw)=\{\eta:\gr^{-1}\Delta(\gr^{*}\mathbb P^\Gw[\eta])
\in L_{P^{*}}(\Gw;\gr\,dx)\} $$
with corresponding norm $$
\norm\eta_{N^{P^{*}}}=\norm {\gr^{-1}\Delta(\gr^{*}\mathbb P^\Gw[\eta])}_{L_{P_{\gr}^{*}}}
$$
and the corresponding capacity $C_{N^{P^{*}}}$. The proof of \rprop{1}-$P$, consequence of Young inequality between Orlicz space is valid with modification. {\bf However}, it appears that the 
full characterization of removable sets cannot be adapted without further 
properties of the function $P^{*}$ like the $\Gd_{2}$-condition. Some results in this directions have been obtained in \cite{Kuz} where a necessary and sufficient condition for removability of boundary set is given, under a very restrictive growth condition on $P$ which reduces the nonlinearity to power-like with limited exponent.

\mysection {Internal measures} 
Several above techniques can be extended to
the following types of problem in which $\mu\in \mathfrak M^b_{+}(\Gw)$:
\begin {equation}
\label {exp1I}\BA {ll} -\Gd u+e^u-1=\gm\qquad&\mbox {in }\;\Gw\\
\phantom{-\Gd +e^u-1}
u=0\quad&\mbox {on }\;\prt\Gw.
\EA
\end {equation}

For this specific problem many interesting results can be found in \cite {BLOP} where the analysis of $\gm$ is made by comparison with the  Hausdorff measure in dimension (N-2), $\CH^{N-2}$. It is proved in particular that if a measure $\gm$ satisfies $\gm\leq 4\gp\CH^{N-2}$, then problem $(\ref{exp2})$ admits a solution, while if $\gm$ charges some Borel set $A$ with Hausdorff dimension less than $N-2$, no solution exists. The results we provide are different and in the Orlicz capacities framework.\medskip

We  define the classes $M_P(\Gw)$ and $M_{P^*}(\Gw)$  similarly to $M_P(\Gw;\gr dx)$ and $M_{P^*}(\Gw;\gr dx)$ except that the measure $\gr dx$ is replaced by the Lebesgue measure $dx$. The Orlicz spaces $L_P(\Gw)$ and $L_{P^*}(\Gw)$ are defined from $M_P(\Gw)$ and $M_{P^*}(\Gw)$ and endowed with the respective Luxemburg norms $\norm {\;}_{P}$ and $\norm {\;}_{P^*}$. We put
\begin {equation}\label{HL0}
\Gd^{L\ln L}(\Gw):=\{\eta\in W^{1,1}_0(\Gw):\Gd\eta\in L_{P^*}(\Gw)\},
\end{equation}
with natural norm
\begin {equation}\label{HL0'}
\norm\eta_{\Gd^{L\ln L}}:=\norm{\eta}_{L^1}+\norm{\Gd\eta}_{L_{P^*}}.
\end{equation}
The norm in  $M_{P^*}(\Gw)$ can be characterized using the Hardy-Littlewood maximal function $f\mapsto M_{Q_0}[f]$ since
\begin {equation}\label{HL1}
\norm f_{L\ln L}:= \int_{Q_{0}}M_{Q_{0}}[f](x) dx\approx\norm f_{L_{P^{*}}}.
\end{equation}
Since $P^*$ satisfies the $\Gd_2$-condition, $C^\infty_0(\Gw)$ is dense in $\Gd^{L\ln L}(\Gw)$. Inequality $(\ref {dual})$ becomes
\begin {equation}
\label {dualI}
\abs{\int_{\Gw}\eta d\gm}=\abs{\int_{\Gw}\eta \Gd \mathbb G^\Gw[\gm]\,dx}
=\abs{\int_{\Gw}\mathbb G^\Gw[\gm]\Gd \eta \,dx}
\leq \norm{\mathbb G^\Gw[\gm]}_{L_{P}}
\norm{\Gd\eta}_{L_{P^{*}}},
\end{equation}
for $\eta\in C^{1,1}_{c}(\bar\Gw)$. We define the $C_{\Gd^{L\ln L}}$-capacity of a compact
subset $K$ of 
$\prt\Gw$ by
\begin {equation}
\label {GdLlnL3} C_{\Gd^{L\ln L}}(K)=\inf\{\norm {\Gd\eta}_{L_{P^{*}}}:\eta\in
C_{c}^{2}(\Gw), 
\eta\geq 0 ,\eta\geq 1\mbox { in a neighborhood of }K\},
\end{equation}
By the min-max theorem there holds
\begin {equation}
\label {GdLlnL3'} C_{\Gd^{L\ln L}}(K)=\sup\{\gm (K):\gm\in \mathfrak M^b_+(\Gw), \gm(K^c)=0,\norm{\BBG^\Gw[\gm]}\leq 1\}.
\end{equation}

\nind\Remark The characterization of the $C_{\Gd^{L\ln L}}$-capacity is not simple, however, by a result of \cite [Th1]{CS}, there holds
\begin {equation}
\label {wL}
\norm {D^2\eta}_{L^{1,\infty}}\leq C\norm{\Gd\eta}_{L\ln L}\qquad\forall \eta\in C^{1,1}_c(\overline\Gw)
\end{equation}
where $L^{1,\infty}(\Gw)$ denotes the weak $L^1$-space, that is the space of all measurable functions $f$ defined in $\Gw$ satisfying
 \begin {equation}
\label {wL1}
meas\left(\{x\in\Gw:|f(x)|>t\}\right)\leq \myfrac{c}{t},
\qquad\forall t>0
\end{equation}
and $\norm {f}_{L^{1,\infty}}$ is the smallest constant such that $(\ref{wL1})$ holds.

\bdef{intm} The space of all bounded measures in $\Gw$ such that $\mathbb G^\Gw[\gm]\in L_{P}(\Gw) $
is denoted by $B^{\exp}(\Gw)$, with norm 
\begin {equation}
\label {Bexp}
\norm {\gm}_{B^{\exp}}=\norm {\mathbb G^\Gw[\gm]}_{L_{P}}.
\end{equation}
The subset of measures such that $\exp(\mathbb G^\Gw[\gm])\in L^{1}(\Gw) $
is denoted by $\mathfrak M^{\exp}(\Gw)$.
\es

Thus \rprop{1} and Theorem B are valid under  the form

\bprop{1'} If $\gm\in B_{+}^{\exp}(\Gw) $, it does 
not charge Borel subsets with $C_{\Gd^{L\ln L}}$-capacity zero.\es

\bth{B'} Let $\gm\in\mathfrak M_+(\Gw)$ be a good measure, then $\gm$ vanishes on Borel subset $E$ with zero $C_{\Gd^{L\ln L}}$-capacity.
\es

Theorem C has the following counter part\medskip

\bth{C'} Let $K\subset\Gw$ be compact. Any solution of  
\begin {equation}
\label {expX}\BA {ll} -\Gd u+e^u-1=0\qquad&\mbox {in }\;\Gw\setminus K\\
\phantom{-\Gd +e^u-1}
u=0\quad&\mbox {on }\;\prt\Gw,
\EA
\end {equation}
vanishes identically in $\Gw$ if and only if $C_{\Gd^{L\ln L}}(K)=0$.
\es
\Proof Let $u\in C(\Gw\setminus K)$ be a solution of $(\ref {exp1})$ which is zero on
$\prt\Gw$. Let $\{\eta_n\}\subset C^2(\Gw)$ such that $0\leq \eta_n\leq1$, $\eta_n=1$ in a
neighborhood $\CV$ of $K$ and $(\ref {AM1})$ holds. Put $\gr_K(x)=\dist (x,K)$. Then, as a
consequence of Keller-Osserman 
estimate for this type of nonlinearity (see \cite {VV}), there holds 
$$
u(x)\leq C\ln (2/\gr_K (x))+D. $$
Put $\gth_n =1-
\eta_n$. Then the function $\gz_n=\phi_{1}\theta_n$ ($\phi_{1}$ being the 
first eigenfunction of $-\Gd$) is an admissible test function for $u$, and $$
\int_\Gw\left(-u\Gd\gz_n+(e^u-1)\gz_n\right) dx=0. $$
We derive
\begin {eqnarray*}
-\int_\Gw u\Gd\gz_n\,dx=-\int_\Gw \gz_n^{-1}\Gd\gz_n\,u\,dx
\phantom {-\int_\Gw u\Gd\gz_n\,dx-------------------,,}\\
\geq -2^{- 1}\int_\Gw(e^u-1-u)\,dx-\int_\Gw Q(\Gd(\gz_n)\,dx.
\phantom {---------}
\end {eqnarray*}
Therefore
\begin {eqnarray}
\label {est}
\int_\Gw(e^u-1-u)\gz_n\,dx
\leq 2C\int_\Gw \abs{\Gd\gz_n}\ln (1+\abs{\Gd\gz_n})dx,
\end {eqnarray}
Since the right-hand side goes to zero when $n\to\infty$, the conclusion follows.\qeda

\subsection{More on good measures}
The main characterization of good measures is the following
\bth{gm} Assume $\gm$ is a positive good measure, then there exists an increasing sequence $\{\gm_n\}\subset B_+^{exp}(\Gw)$ which converges weakly to $\gm$.\es

The proof will necessitate several intermediate results which are classical in the framework of Lebesgue measure or Bessel capacities, but appear to be new for Orlicz capacities.\medskip

 \blemma {L1}Let $K\subset\Gw$, then $C_{\Gd^{L\ln L}}(K)=0$ if and only if there exists $\eta\in \Gd^{L\ln L}(\Gw)$ such that 
 $\eta\geq 0$ and $K\subset\{y\in\Gw:\eta(y)=\infty\}$.
 \es
\Proof By the definition of the capacity, for any $\gl>0$ and $\eta\in \Gd^{L\ln L}(\Gw)$, $\eta\geq 0$,  
\begin {equation}\label {check1} C_{\Gd^{L\ln L}}\left(\{y\in\Gw:\eta(y)\geq\gl\}\right)\leq \myfrac{1}{\gl}
\norm{\eta}_{\Gd^{L\ln L}}.
\end{equation}
This implies
$$C_{\Gd^{L\ln L}}\left(\{y\in\Gw:\eta(y)=\infty\}\right)=0.
$$
\qeda
 \blemma  {L2} Suppose $\{\eta_j\}$ is a Cauchy sequence in $\Gd^{L\ln L}(\Gw)$. Then there exist a subsequence $\{\eta_{j_\ell}\}$ and $\eta\in \Gd^{L\ln L}(\Gw)$ such that 
 $$\lim_{i_\ell\to\infty}\eta_{j_\ell}=\eta,$$
 uniformly outside an open subset of arbitrary small $C_{\Gd^{L\ln L}}$-capacity.
 \es
 \Proof By \rlemma {L1}, $\eta_{j}$ and $\eta$ are finite outside a set $F$ with zero $C_{\Gd^{L\ln L}}$-capacity.  There exists a subsequence $\{\eta_{j_\ell}\}$ such that 
 $$\norm{\eta_{j_\ell}-\eta}_{\Gd^{L\ln L}}\leq 2^{-2\ell}.
 $$
 Put $E_\ell=\{y\in\Gw:\eta_{j_\ell}-\eta(y)\geq 2^{-\ell}\}$. By $(\ref{check1})$
 $C_{\Gd^{L\ln L}}\left(E_\ell\right)\leq  2^{-\ell}$, and if $G_m=\cup_{\ell\geq m}E_\ell$, there holds $C_{\Gd^{L\ln L}}(G_m)\leq 2^{1-m}$. Therefore
 $$C_{\Gd^{L\ln L}}(\cap_{ m\geq 1}G_m)=0.$$
Since for any $y\notin G_m\cup F$, there holds
 $$\abs{(\eta_{j_\ell}-\eta)(y)}\leq 2^{-\ell},
 $$
the claim follows.\qeda
 \medskip

 \blemma  {L3}  If $\eta\in \Gd^{L\ln L}(\prt\Gw)$ it has a unique quasi-continuous representative with respect to the capacity $C_{\Gd^{L\ln L}}$.
 \es
 \Proof Uniqueness is clear as in the Bessel capacity case \cite[Chap 6]{AdHe}. Let $\{\eta_j\}\subset C_0^2(\Gw)$ be a sequence which converges to $\eta$ in $\Gd^{L\ln L}(\Gw)$. Then there exists a subsequence $\{\eta_{j_\ell}\}$ such that $\eta_{j_\ell}$  converges
 to $\eta$ uniformly on the complement of an open set of arbitrarily small  $C_{\Gd^{L\ln L}}$-capacity. This is the claim.
 \qeda\medskip 

\nind{\it Proof of \rth{gm}.} The method is adapted from  \cite[Th 8]{FeDel}, \cite[Lemma 4.2]{BaPi}. By \rlemma{L3} we can define the functional $h$ on $\Gd^{L\ln L}(\Gw)$ by
$$h(\eta)
=\int_{\Gw}\overline\eta_+ d\gm\qquad\forall\eta\in \Gd^{L\ln L}(\Gw),
$$
where $\overline {\eta}$ stands for the $C_{\Gd^{L\ln L}}$-quasi-continuous representative of $\eta$. Notice that we can write
$$\BA {l}h(\eta)=-\myint{\Gw}{}\Gd\BBG^\Gw[\gm]\eta dx=-\myint{\Gw}{}\BBG^\Gw[\gm]\Gd\eta dx
\EA$$

The following steps are similar to the previous proofs:\smallskip

\nind {\it Step 1-} The functional $h$ is convex, positively homogeneous and l.s.c. The convexity and the homogeneity are clear. If $\eta_n\to\eta$ in $\Gd^{L\ln L}(\prt\Gw)$, then by \rlemma{L3} we can extract a subsequence which is converging everywhere except for a set with zero capacity. The conclusion follows from Fatou's lemma.\smallskip

\nind {\it Step 2-} Since $L_{P}(\Gw)$ is the dual space of $L_{P^*}(\Gw)$, for any continuous linear form $\ga$ on $\Gd^{L\ln L}(\Gw)$ there exists $\gb\in L_{P}(\Gw)$ such that
$$\ga(\eta)=-\myint{\Gw}{}\gb\Gd\eta dx\qquad\forall\eta\in \Gd^{L\ln L}(\Gw).
$$
Therefore, in the sense of distributions there holds
$$\ga(\eta)=-\langle\Gd\gb,\eta\rangle\qquad\forall\eta\in C^\infty_0(\Gw).
$$

\nind {\it Step 3-} By the geometric Hahn-Banch theorem, $h$ is the upper convex hull of the continuous linear functionals on  $\Gd^{L\ln L}(\prt\Gw)$ it dominates. Fix a function $\eta_0\in C^\infty_0(\Gw)$ and $\ge>0$, there exists a continuous linear form
$\ga$ on  $\Gd^{L\ln L}(\Gw)$ and constants $a,b$ such that 
$$a+bt+\ga(\eta)\leq 0\qquad\forall (\eta,t)\in \CE:=\{(\eta,t)\in \Gd^{L\ln L}(\Gw)\ti\BBR:h(\eta)\leq t\},
$$
and
$$a+b(h(\eta_0)-\ge)+\ga(\eta_0)> 0.
$$
The same ideas as in \cite[Lemma 4.2]{BaPi} yields successively to $a=0$ and $b<0$. If we put $\gs(\eta)=-b^{-1}\ga(\eta)$ we derive $\gs(\eta)\leq h(\eta)$ for all $\eta\in \Gd^{L\ln L}(\Gw)$. This implies in particular that $\gs(\eta)\leq 0$ if $\eta\leq 0$, thus
$\gs$ is a positive linear form on $\Gd^{L\ln L}(\Gw)$. Therefore there exists a Radon measure $\gn$ on $\Gw$ and $\gb\in L_P(\Gw)$ such that $-\Gd\gb=\gn$, 
$0\leq \gn\leq \gm$ and 
$$\myint{\Gw}{}\eta_0 d\gm\leq \ge+\myint{\Gw}{}\eta_0 d\gn.
$$

\nind {\it Step 4-} Considering an increasing sequence of compact sets $K_j$ such that $K_j\subset \overset{o}{K}_{j+1}$ and $\cup_jK_j=\Gw$, we construct for each $j\in\BBN^*$ a Radon measure $\gn_j$ and $\gb_j\in L_P(\Gw)$ such that $-\Gd\gb_j=\gn_j$, 
$0\leq \gn_j\leq \gm$ and 
$$\myint{K_j}{}d\gm\leq j^{-1}+\myint{K_j}{}d\gn_j.
$$
At last we can assume that the sequence $ \{\gn_j\}$ is increasing since if $-\Gd\gb_j=\gn_j$ for $j=1,2$, then 
$$-\Gd\gb_{1,2}=\sup \{\gn_1,\gn_2\}\leq \gn_1+\gn_2=-\Gd\gb_1-\Gd\gb_2
$$
thus $\gb_{1,2}\in L_P(\Gw)$. Iterating this process, we can replace the sequence $\{\gn_j\}$ by $ \{\gn'_j\}:=\{\gn_1,\sup \{\gn_1,\gn_2\},\sup\{\gn_3,\sup \{\gn_1,\gn_2\}\}...\}$. The sequence $ \{\gn'_j\}$ is increasing, converges to $\gm$ and since $\gb_j=\BBG^{\Gw}[\gn'_j]$ with $\gb_j\in L_P(\Gw)$, $\gn'_j$ belongs to $B^{exp}(\Gw)$.\qeda\medskip

As a consequence of this result and the characterization of linear functionals over $L\ln L(\Gw)$, the following result holds.

\bcor{end} Assume $\gm$ is a bounded positive good measure in $\Gw$, then ther exists an increasing sequence of positive measures  $\gn_j$ in $\Gw$ and positive real numbers $\gth_j$ such that $\gn_j$ to $\gm$ in the weak sense of measures and $\exp{\left(\gth_j\BBG^{\Gw}[\gn_j]\right)}\in L^1(\Gw)$.
\es

\end {document}